\documentclass{ article}

\begin{document}
\vspace*{.5cm}
\begin{center}
{\Large{\bf   Slant submersions from almost Hermitian manifolds}}\\
\vspace{.5cm}
 { Bayram \d{S}ahin} \\
\end{center}

\vspace{.5cm}
\begin{center}
{\it Inonu University, Department of Mathematics, 44280,
Malatya-Turkey. E-mail:bsahin@inonu.edu.tr}
\end{center}
\vspace{.5cm}

\noindent {\bf Abstract.} {\small We introduce slant submersions
from almost Hermitian manifolds onto Riemannian manifolds. We give
examples, investigate the geometry of foliations which are arisen
from the definition of a Riemannian submersion and check the
harmonicity of such submersions. We also find necessary and
sufficient conditions for a slant submersion to be totally geodesic.
Moreover, we obtain a decomposition theorem for the total manifold
of
such submersions.}\\



\pagestyle{myheadings}

\section*{1.Introduction}

  \setcounter{equation}{0}
\renewcommand{\theequation}{1.\arabic{equation}}
\markboth{Slant submersions}{Slant
submersions}{\thispagestyle{plain}}

Let $\bar{M}$ be a Kaehler manifold with complex structure $J$ and
$M$ is a Riemannian manifold isometrically immersed in $\bar{M}.$ We
note that  submanifolds of a Kaehler manifold determined by the
behavior of the tangent bundle of the submanifold under the action
of the complex structure of the ambient manifold. A submanifold $M$
is called holomorphic (complex) if $J (T_p M)\subset T_p M$, for
every $p \in M$, where $T_p M$  denotes the tangent space to $M$ at
the point $p .$ $M$ is called totally real if $J(T_p M) \subset T_p
M^{\perp}$ for every $p \in M,$ where $T_p M^{\perp}$ denotes the
normal space to $M$ at the point $p$.    As a generalization of
holomorphic and totally real submanifolds, slant submanifolds were
introduced by Chen in \cite{CB}. We recall that the submanifold $M$
is called slant \cite{CB} if for all
    non-zero vector $X$ tangent to $M$ the angle $\theta (X)$
    between $J X$ and $T_p M$ is a constant, i.e, it does not depend
    on the choice of $p \in M$ and $X \in T_p M.$.\\

On the other hand, Riemannian submersions between Riemannian
manifolds were studied by O'Neill \cite{O'Neill} and Gray
\cite{Gray}. Later such submersions were considered between
manifolds with differentiable structures. As an  analogue of
holomorphic submanifolds, Watson defined almost Hermitian
submersions between almost Hermitian manifolds and he showed that
the base manifold and each fiber have the same kind of structure as
the total space, in most cases \cite{Watson}. We note that almost
Hermitian submersions have been extended to the almost contact
manifolds \cite{Domingo}, locally conformal K\"{a}hler
manifolds\cite{lck}
and quaternion K\"{a}hler manifolds \cite{Ianus}.\\

Let $M$ be a complex $m-$dimensional almost Hermitian manifold with
Hermitian metric $g_M$ and almost complex structure $J_M$ and $N$ be
a complex $n-$dimensional almost Hermitian manifold with Hermitian
metric $g_N$ and almost complex structure $J_N$. A Riemannian
submersion $F:M\longrightarrow N$ is called an almost Hermitian
submersion if $F$ is an almost complex mapping, i.e.,
$F_*J_M=J_NF_*$. The main result of this notion is that the vertical
and horizontal distributions are $J_M-$ invariant. On the other
hand, Riemannian submersions from almost Hermitian manifolds onto
Riemannian manifolds have been studied by many authors under the
assumption that the vertical spaces of such submersions are
invariant with respect to the complex structure. For instance,
Escobales\cite{Escobales} studied Riemannian submersions from
complex projective space onto a Riemannian manifold under the
assumption that the fibers are connected, complex, totally geodesic
submanifolds. One can see that this assumption implies that the
vertical distribution is invariant.\\

Recently, we introduce anti-invariant Riemannian submersions from
almost Hermitian manifolds onto Riemannianian manifolds and
investigate the geometry of such submersions \cite{Sahin}. In this
paper, as a generalization of Hermitian submersions and
anti-invariant submersions, we define and study slant submersions
from almost Hermitian manifolds onto
Riemannian manifolds.\\

The paper is organized as follows: In section~2, we present the
basic information needed for this paper. In section~3, we give
definition of slant Riemannian submersions, provide examples and
investigate the geometry of leaves of the distributions. We also
obtain necessary and sufficient conditions for such submersions to
be totally geodesic. Moreover we give a necessary condition for
slant submersions to be harmonic and obtain a decomposition theorem.

\section*{2.Preliminaries}
  \setcounter{equation}{0}
\renewcommand{\theequation}{2.\arabic{equation}}
In this section, we define almost Hermitian manifolds, recall the
notion of Riemannian submersions between Riemannian manifolds and
give a brief review
of basic facts of Riemannian submersions.\\

Let ($\bar{M}, g$) be an almost Hermitian manifold. This means
\cite{Yano-Kon} that $\bar{M}$ admits a tensor field $J$ of type (1,
1) on $\bar{M}$ such that, $\forall X, Y \in \Gamma(T\bar{M})$, we
have
\begin{equation}
J^2=-I, \quad g(X, Y)=g(JX, JY) \label{eq:2.1}.
\end{equation}
 An almost Hermitian manifold $\bar{M}$ is called  K\"{a}hler manifold if
\begin{equation}
(\bar{\nabla}_XJ)Y=0,\forall X,Y \in \Gamma(T\bar{M}),
\label{eq:2.2}
\end{equation}
where $\bar{\nabla}$ is the Levi-Civita connection on $\bar{M}. $

Let $(M^m,g_{_M})$ and $(N^n,g_{_N})$ Riemannian manifolds, where
$dim(M)=m$, $dim(N)=n$ and $m>n$. A Riemannian submersion
$F:M\longrightarrow N$ is a map of $M$ onto $N$ satisfying the
following axioms:
\begin{itemize}
    \item [(S1)] $F$ has maximal rank.
    \item [(S2)] The differential $F_*$ preserves the lenghts of
    horizontal vectors.
\end{itemize}
For each $q\in N$, $F^{-1}(q)$ is an $(m-n)$ dimensional submanifold
of $M$.The submanifolds $F^{-1}(q)$, $q\in N$, are called fibers. A
vector field on $M$ is called vertical if it is always tangent to
fibers. A vector field on $M$ is called horizontal if it is always
orthogonal to fibers. A vector field $X$ on $M$ is called basic if
$X$ is horizontal and $F-$ related to a vector field $X_*$ on $N$,
i.e., $F_*X_p=X_{*F(p)}$ for all $p \in M$. Note that we denote the
projection morphisms on the distributions $ker F_*$ and $(kerf
F_*)^\perp$ by $\mathcal{V}$ and $\mathcal{H}$, respectively.\\

We recall the following lemma from O'Neill \cite{O'Neill}.\\
\noindent{\it{\bf Lemma~2.1~}}{\it Let $F:M \longrightarrow N$ be a
Riemannian submersion between Riemannian manifolds and $X, Y$ be
basic vector fields of $M$. Then }
\begin{itemize}
    \item [(a)] $g_{_M}(X, Y)=g_{_N}(X_*,Y_*)\circ F$
    \item [(b)] {\it the horizontal part $[X,Y]^{\mathcal{H}}$ of
    $[X,Y]$ is a basic vector field and corresponds to
    $[X_*,Y_*]$,i.e., $F_*([X,Y]^{\mathcal{H}})=[X_*,Y_*]$}.
    \item [(c)] $[V,X]$ {\it is vertical for any  vector field $V$ of
    $ker F_*$.}
    \item [(d)] $(\nabla^{^M}_XY)^{\mathcal{H}}$ {\it is the basic vector field  corresponding to $\nabla^{^N}_{X_*}Y_*$.}
\end{itemize}
The geometry of Riemannian submersions is characterized by O'Neill's
tensors $\mathcal{T}$ and $\mathcal{A}$ defined for vector fields
$E, F$ on $M$ by
\begin{equation}
\mathcal{A}_E
F=\mathcal{H}\nabla_{\mathcal{H}E}\mathcal{V}F+\mathcal{V}\nabla_{\mathcal{H}E}\mathcal{H}F
\label{eq:2.3}
\end{equation}
\begin{equation}
\mathcal{T}_E
F=\mathcal{H}\nabla_{\mathcal{V}E}\mathcal{V}F+\mathcal{V}\nabla_{\mathcal{V}E}\mathcal{H}F,
\label{eq:2.4}
\end{equation}
where $\nabla$ is the Levi-Civita connection of $g_{_M}$. It is easy
to see that a Riemannian submersion $F:M\longrightarrow N$ has
totally geodesic fibers if and only if $\mathcal{T}$ vanishes
identically. For any $E \in \Gamma(TM)$, $\mathcal{T}_E$ and
$\mathcal{A}_E$ are skew-symmetric operators on $(\Gamma(TM),g)$
reversing the horizontal and the vertical distributions. It is also
easy to see that $\mathcal{T}$ is vertical,
$\mathcal{T}_E=\mathcal{T}_{\mathcal{V}E}$ and $\mathcal{A}$ is
horizontal, $\mathcal{A}=\mathcal{A}_{\mathcal{H}E}$. We note that
the tensor fields $\mathcal{T}$ and $\mathcal{A}$ satisfy
\begin{eqnarray}
\mathcal{T}_UW&=&\mathcal{T}_WU, \forall U,W \in \Gamma(ker F_*)\label{eq:2.5}\\
\mathcal{A}_XY&=&-\mathcal{A}_YX=\frac{1}{2}\mathcal{V}[X,Y],
\forall X,Y\in \Gamma((ker F_*)^{\perp}).\label{eq:2.6}
\end{eqnarray}
On the other hand, from (\ref{eq:2.3}) and (\ref{eq:2.4}) we have
\begin{eqnarray}
\nabla_VW&=&\mathcal{T}_VW+\hat{\nabla}_V W \label{eq:2.7}\\
\nabla_VX&=&\mathcal{H}\nabla_VX+\mathcal{T}_VX \label{eq:2.8}\\
\nabla_XV&=&\mathcal{A}_XV+\mathcal{V}\nabla_XV \label{eq:2.9}\\
\nabla_XY&=&\mathcal{H}\nabla_XY+\mathcal{A}_XY \label{eq:2.10}
\end{eqnarray}
for $X,Y \in \Gamma((ker F_*)^{\perp})$and $V,W \in \Gamma(ker
F_*)$, where $\hat{\nabla}_VW=\mathcal{V}\nabla_VW$. If $X$ is
basic, then $\mathcal{H}\nabla_VX=\mathcal{A}_XV$.

Finally, we recall the notion of harmonic maps between Riemannian
manifolds. Let $(M, g_{_M})$ and $(N, g_{_N})$ be Riemannian
manifolds and suppose that $\varphi: M\longrightarrow N$ is a smooth
map between them. Then the differential $\varphi_{*}$ of $\varphi$
can be viewed a section of the bundle $Hom(TM,
\varphi^{-1}TN)\longrightarrow M,$ where $\varphi^{-1}TN$ is the
pullback bundle which has fibers $(\varphi^{-1}TN)_p=T_{\varphi(p)}
N, p \in M.$ $Hom(TM, \varphi^{-1}TN)$ has a connection $\nabla$
induced from the Levi-Civita connection $\nabla^M$ and the pullback
connection. Then the second fundamental form of $\varphi$ is given
by
\begin{equation}
(\nabla \varphi_{*})(X, Y)=\nabla^{\varphi}_X
\varphi_{*}(Y)-\varphi_{*}(\nabla^M_X Y) \label{eq:2.11}
\end{equation}
for $X, Y \in \Gamma(TM)$, where $\nabla^{\varphi}$ is the pullback
connection. It is known that the second fundamental form is
symmetric. A smooth map $\varphi: (M, g_{_M}) \longrightarrow (N,
g_{_N})$ is said to be harmonic if $trace (\nabla \varphi_{*})=0.$
On the other hand, the tension field of $\varphi$ is the section
$\tau(\varphi)$ of $\Gamma(\varphi^{-1}TN)$ defined by
\begin{equation}
\tau(\varphi)=div\varphi_{*}=\sum^{m}_{i=1} (\nabla
\varphi_{*})(e_i, e_i), \label{eq:2.12}
\end{equation}
where $\{e_1,...,e_m\}$ is the orthonormal frame on $M.$ Then it
follows that $\varphi$ is harmonic if and only if
$\tau(\varphi)=0$, for details, see \cite{Baird-Wood}\\

\section*{3.Slant Submersions}
  \setcounter{equation}{0}
\renewcommand{\theequation}{3.\arabic{equation}}

In this section, we define slant submersions from an almost
Hermitian manifold onto a Riemannian manifold by using the
definition of a slant distribution given in \cite{Carriazo2} ,
investigate the integrability of distributions and obtain a
necessary and sufficient condition for such submersions to be
totally geodesic map. We also investigate the harmonicity of a
slant submersions and obtain a decomposition theorem for the total manifold.\\

\noindent{\bf Definition~3.1.~}{\it Let $F$ be a Riemannian
submersion from an almost Hermitian manifold $(M_1,g_1,J_1)$  onto a
Riemannian manifold $(M_2,g_2)$. If for any non-zero vector $X \in
\Gamma(ker F_*)$, the angle $\theta(X)$ between $JX$ and the space
$ker F_*$ is a constant, i.e. it is independent of the choice of the
point $p \in M_1$ and choice of the tangent vector $X$ in $ker F_*$,
then we say that $F$ is a slant submersion. In this case, the angle $\theta$ is called the slant angle of the slant submersion.}\\

It is known that the distribution $ker F_*$ is integrable. In fact,
its leaves are $F^{-1}(q)$, $q\in M_1$, i.e., fibers. Thus it
follows from above definition that the fibers are slant submanifolds
of $M_1$, for
slant submanifold, \cite{Chen}.\\

We first give some examples of slant submersions.\\

\noindent{\bf Example~1.~} Every Hermitian submersion from an almost
Hermitian manifold onto an almost Hermitian manifold is a slant
submersion with $\theta=0$.\\

\noindent{\bf Example~2.~} Every anti-invariant Riemannian
submersion from an almost Hermitian manifold to a Riemannian
manifold is a slant submersion with $\theta=\frac{\pi}{2}$.\\

A slant submersion is said to be proper if it is neither Hermitian
nor anti-invariant Riemannian submersion.\\

 \noindent{\bf Example~3.~} Consider the following Riemannian
submersion given by
$$
\begin{array}{cccc}
  F: & R^4             & \longrightarrow & R^2\\
     & (x_1,x_2,x_3,x_4) &             & (x_1\,\sin \alpha -x_3\,\cos \alpha, x_4).
\end{array}
$$
Then for any $0< \alpha <\frac{\pi}{2}$, $F$ is a slant submersion
with slant angle $\alpha$.\\

\noindent{\bf Example~4.~}The following Riemannian submersion
defined by
$$
\begin{array}{cccc}
  F: & R^4             & \longrightarrow & R^2\\
     & (x_1,x_2,x_3,x_4) &             & (\frac{x_1 - x_4}{\sqrt{2}},
     x_2)
\end{array}
$$
is a slant submersion with slant angle $\theta=\frac{\pi}{4}$.\\

Let $F$ be a  Riemannian submersion from an almost Hermitian
manifold $(M_1,g_1,J)$ onto a Riemannian manifold $(M_2,g_2)$. Then
for $X \in \Gamma(ker F_*)$, we write
\begin{equation}
JX=\phi X+\omega X, \label{eq:3.1}
\end{equation}
where $\phi X$ and $ \omega X$ are vertical and horizontal parts of
$JX$. Also for $V \in \Gamma((ker F_*)^\perp)$, we have
\begin{equation}
JZ=\mathcal{B}Z+\mathcal{C}Z, \label{eq:3.2}
\end{equation}
where $\mathcal{B}Z$ and $\mathcal{C}Z$ are vertical and horizontal
components of $JZ$. Using (\ref{eq:2.7}), (\ref{eq:2.8}),
(\ref{eq:3.1}) and (\ref{eq:3.3}) we obtain
\begin{eqnarray}
(\nabla_X \omega)Y&=&\mathcal{C}\mathcal{T}_XY-\mathcal{T}_X \phi Y
\label{eq:3.3}\\
(\nabla_X \phi)Y&=&\mathcal{B}\mathcal{T}_XY-\mathcal{T}_X \omega Y,
\label{eq:3.4}
\end{eqnarray}
where
\begin{eqnarray}
(\nabla_X \omega)Y&=&\mathcal{H}\nabla_X \omega Y-\omega
\hat{\nabla}_X Y\nonumber\\
(\nabla_X \phi)Y&=& \hat{\nabla}_X \phi Y -\phi
\hat{\nabla}_XY.\nonumber
\end{eqnarray}
for $X, Y \in \Gamma(ker F_*)$. Let $F$ be a slant submersion from
an almost Hermitian manifold onto a Riemannian manifold, then we say
that $\omega$ is parallel if
$\nabla \omega =0$.\\

The proof of the following result is exactly same with slant
immersions (see \cite{Chen} or \cite{Carriazo} for Sasakian case),
therefore we omit its proof.\\

\noindent{\bf Theorem~3.1.~}{\it Let $F$ be a Riemannian submersion
from an almost Hermitian manifold $(M_1,g_1,J)$ onto a Riemannian
manifold $(M_2,g_2)$. Then $F$ is a slant submersion if and only if
for,  $\phi^2X=\lambda X$, $\lambda \in [-1,0]$ and $X \in
\Gamma(ker F_*)$. If $F$ is a slant submersion, then
$\lambda=-\cos^2 \theta$.}

By using above theorem, it is easy to see that
\begin{eqnarray}
g_1(\phi X, \phi Y)&=&\cos^2 \theta g_1(X,Y) \label{eq:3.5}\\
g_1(\omega X, \omega Y)&=&\sin^2 \theta g_1(X,Y) \label{eq:3.6}
\end{eqnarray}
for any $X, Y \in \Gamma(ker F_*)$. Also by using (\ref{eq:3.5}) and
(\ref{eq:3.6}) we can easily conclude that
$$\{e_1,\sec \theta \phi e_1, e_2,\sec \theta \phi e_2,...,e_n,\sec \theta \phi
e_n\}$$ is an orthonormal frame for $\Gamma(ker F_*)$ and
$$\{\csc \theta \omega e_1, \csc \theta \omega e_2,...,\csc \theta \omega
e_n\}$$ is an orthonormal frame for $\Gamma(J(ker F_*))$. As in
slant immersions, we call such frame an adapted frame for slant
submersions.\\

\noindent{\bf Lemma 3.1.~}{\it Let $F$ be a slant submersion from a
K\"{a}hler manifold onto a Riemannian manifold. If $\omega$ is
parallel then we have}
\begin{equation}
\mathcal{T}_{\phi X}\phi X=-\cos^2 \theta \mathcal{T}_X X
\label{eq:3.7}
\end{equation}
{\it for $X \in \Gamma(ker F_*)$.}\\

\noindent{\bf Proof.~} If $\omega$ is parallel, then from
(\ref{eq:3.3}) we have $\mathcal{C}\mathcal{T}_X Y=\mathcal{T}_X
\phi Y$ for $X, Y \in \Gamma(ker F_*)$. Interchanging the role of
$X$ and $Y$, we get $\mathcal{C}\mathcal{T}_Y X=\mathcal{T}_Y \phi
X$. Thus we have
$$\mathcal{C}\mathcal{T}_X Y-\mathcal{C}\mathcal{T}_Y X=\mathcal{T}_X
\phi Y-\mathcal{T}_Y \phi X.$$ Using (\ref{eq:2.5}) we derive
\begin{equation}
\mathcal{T}_X \phi Y=\mathcal{T}_Y \phi X. \label{eq:3.8}
\end{equation}
Then substituting $Y$ by $\phi X$ we get $\mathcal{T}_X \phi^2
X=\mathcal{T}_{\phi X} \phi X.$ Finally using Theorem 3.1 we obtain
(\ref{eq:3.7}).\\

We now give a sufficient condition for a slant submersion to be harmonic.\\

\noindent{\bf Theorem~3.2.~} {\it Let $F$ be a slant submersion from a K\"{a}hler manifold onto a Riemannian manifold. If $\omega$ is parallel then $F$ is a harmonic map.}\\

\noindent{\bf Proof.~} Since
\begin{equation}
(\nabla F_*)(Z_1,Z_2)=0 \label{eq:3.9}
\end{equation}
for $Z_1,Z_2 \in \Gamma((ker F_*)^\perp)$. A slant submersion $F$ is
harmonic if and only if $\sum^{n}_{i=1}(\nabla
F_*)(\tilde{e}_i,\tilde{e}_i)=-\sum^{n}_{i=1}F_*(\mathcal{T}_{\tilde{e}_i}\tilde{e}_i)=0$,
where $\{ \tilde{e}_i\}^n_{i=1}$ is an orthonormal basis of $ker
F_*$. Thus using the adapted frame of slant submersions we can write
$$\tau=-\sum_{i=1}^{\frac{n}{2}}F_*(\mathcal{T}_{e_i}e_i+\mathcal{T}_{\sec
\theta \phi e_i}\sec \theta \phi e_i).$$ Hence we have
$$\tau=-(\sum_{i=1}^{\frac{n}{2}}F_*(\mathcal{T}_{e_i}e_i+\sec^2 \theta \mathcal{T}_{\phi e_i}\phi e_i)).$$
Then using (\ref{eq:3.7}) we arrive at
$$\tau=-(\sum_{i=1}^{\frac{n}{2}}F_*(\mathcal{T}_{e_i}e_i- \mathcal{T}_{e_i}e_i))=0$$
which shows that $F$ is harmonic.

We now investigate the geometry of the leaves of the distributions
$ker F_*$ and $(ker F_*)^\perp$.\\

\noindent{\bf Theorem~3.3.~}{\it Let $F$ be a slant submersion from
a K\"{a}hler manifold $(M_1,g_1,J_1)$ onto a Riemannian manifold
$(M_2,g_2)$. Then the distribution $ker F_*$ defines a totally
geodesic foliation on $M_1$ if and only if}
$$g_1(\mathcal{H}\nabla_X \omega \phi Y,Z)=g_1(\mathcal{H}\nabla_X
\omega Y,\mathcal{C}Z)+g_1(\mathcal{T}_X \omega Y,\mathcal{B}Z)$$
{\it for $X, Y \in \Gamma(ker F_*)$ and $Z \in \Gamma((ker F_*)^\perp$.}\\

\noindent{\bf Proof.~} For $X, Y \in \Gamma(ker F_*)$ and $Z \in
\Gamma((ker F_*)^\perp$, from (\ref{eq:2.1}) and (\ref{eq:3.1}) we
have

$$g_1(\nabla_X Y,Z)=g_1(\nabla_X \phi Y,JZ)+g_1(\nabla_X \omega
Y,JZ).$$ Using (\ref{eq:2.1}), (\ref{eq:3.1}) and (\ref{eq:3.2}) we
get
\begin{eqnarray}
g_1(\nabla_X Y,Z)&=&-g_1(\nabla_X \phi^2 Y,Z)-g_1(\nabla_X \omega
\phi Y,Z)\nonumber\\
&+&g_1(\nabla_X \omega Y,\mathcal{B}Z)+g_1(\nabla_X \omega
Y,\mathcal{C}Z).\nonumber
\end{eqnarray}
Then from (\ref{eq:2.8}) and Theorem 3.1 we obtain
\begin{eqnarray}
g_1(\nabla_X Y,Z)&=&\cos^2 \theta g_1(\nabla_X
Y,Z)-g_1(\mathcal{H}\nabla_X \omega
\phi Y,Z)\nonumber\\
&+&g_1(\mathcal{T}_X \omega Y,\mathcal{B}Z)+g_1(\mathcal{H}\nabla_X
\omega Y,\mathcal{C}Z).\nonumber
\end{eqnarray}
Hence we have
\begin{eqnarray}
\sin^2 \theta g_1(\nabla_X Y,Z)&=&-g_1(\mathcal{H}\nabla_X \omega
\phi Y,Z)\nonumber\\
&+&g_1(\mathcal{T}_X \omega Y,\mathcal{B}Z)+g_1(\mathcal{H}\nabla_X
\omega Y,\mathcal{C}Z)\nonumber
\end{eqnarray}
which proves assertion.\\

In a similar way we have the following.\\

\noindent{\bf Theorem~3.4.~}{\it Let $F$ be a slant submersion from
a K\"{a}hler manifold $(M_1,g_1,J_1)$ onto a Riemannian manifold
$(M_2,g_2)$. Then the distribution $(ker F_*)^\perp$ defines a
totally geodesic foliation on $M_1$ if and only if}
$$g_1(\mathcal{H}\nabla_{Z_1}Z_2,\omega \phi
X)=g_1(\mathcal{A}_{Z_1}\mathcal{B}Z_2+\mathcal{H}\nabla_{Z_1}\mathcal{C}Z_2,\omega
X)$$
{\it for $X\in \Gamma(ker F_*)$ and $Z_1,Z_2 \in \Gamma((ker F_*)^\perp$.}\\

From Theorem 3.3 and Theorem 3.4 we have the following result.\\

\noindent{\bf Corollary~3.1.~}{\it Let $F$ be a slant submersion
from a K\"{a}hler manifold $(M_1,g_1,J_1)$ onto a Riemannian
manifold $(M_2,g_2)$. Then $M_1$ is locally a product Riemannian
manifold if and only if}
$$g_1(\mathcal{H}\nabla_{Z_1}Z_2,\omega \phi
X)=g_1(\mathcal{A}_{Z_1}\mathcal{B}Z_2+\mathcal{H}\nabla_{Z_1}\mathcal{C}Z_2,\omega
X)$$ and
$$g_1(\mathcal{H}\nabla_X \omega \phi Y,Z_1)=g_1(\mathcal{H}\nabla_X
\omega Y,\mathcal{C}Z_1)+g_1(\mathcal{T}_X \omega
Y,\mathcal{B}Z_1)$$
 {\it for $X, Y,\in \Gamma(ker F_*)$ and $Z_1,Z_2 \in
\Gamma((ker F_*)^\perp$.}\\

Finally we give necessary and sufficient conditions for a slant
submersion to be totally geodesic. Recall that  a differentiable map
$F$ between Riemannian manifolds $(M_1,g_1)$ and $(M_2,g_2)$ is
called a totally geodesic map if $(\nabla
F_*)(X,Y)=0$ for all $X, Y \in \Gamma(TM_1)$.\\

\noindent{\bf Theorem~3.5.~}{\it Let $F$ be a slant submersion from
a K\"{a}hler manifold $(M_1,g_1,J_1)$ onto a Riemannian manifold
$(M_2,g_2)$. Then $F$ is totally geodesic if and only if}
$$g_1(\mathcal{T}_X \omega Y,\mathcal{B}Z_1)+g_1(\mathcal{H}\nabla_X
\omega Y,\mathcal{C}Z_1)=g_1(\mathcal{H}\nabla_X \omega \phi
Y,Z_1)$$ and
$$g_1(\mathcal{A}_{Z_1}BZ_2+\mathcal{H}\nabla_{Z_1}\mathcal{C}Z_2,
\omega X)=-g_1(\mathcal{H}\nabla_{Z_1}\omega \phi X,Z_2)$$ {\it for
$Z_1,Z_2 \in \Gamma((ker F_*)^\perp)$ and $X, Y \in \Gamma(ker
F_*)$.}

\noindent{\bf Proof.~} For $X, Y \in \Gamma(ker F_*)$ and $Z_1 \in
\Gamma((ker F_*)^\perp)$, since $F$ is a Riemannian submersion, from
(\ref{eq:2.1}), (\ref{eq:3.1}) and (\ref{eq:3.2}) we have
$$g_2((\nabla F_*)(X,Y),F_*Z_1)=-g_1(\nabla_X J\phi
Y,Z)+g_1(\nabla_X \omega Y,JZ).$$ Using again (\ref{eq:3.1}) and
(\ref{eq:3.2}) we get
\begin{eqnarray}
g_2((\nabla F_*)(X,Y),F_*Z_1)&=&-g_1(\nabla_X \phi^2
Y,Z)-g_1(\nabla_X\omega \phi Y,Z)\nonumber\\
&+&g_1(\nabla_X \omega Y,\mathcal{B}Z)+g_1(\nabla_X \omega
Y,\mathcal{C}Z).\nonumber
\end{eqnarray}
Then Theorem 3.1, (\ref{eq:2.7}) and (\ref{eq:2.8})imply that
\begin{eqnarray}
g_2((\nabla F_*)(X,Y),F_*Z_1)&=&\cos^2 \theta g_1(\nabla_X \phi^2
Y,Z)-g_1(\mathcal{\nabla}_X\omega \phi Y,Z)\nonumber\\
&+&g_1(\mathcal{T}_X \omega Y,\mathcal{B}Z)+g_1(\mathcal{H}\nabla_X
\omega Y,\mathcal{C}Z).\nonumber
\end{eqnarray}
Hence we obtain
\begin{eqnarray}
\sin^2 \theta g_2((\nabla F_*)(X,Y),F_*Z_1)&=& g_1(\nabla_X \phi^2
Y,Z)-g_1(\mathcal{\nabla}_X\omega \phi Y,Z)\nonumber\\
&+&g_1(\mathcal{T}_X \omega Y,\mathcal{B}Z)+g_1(\mathcal{H}\nabla_X
\omega Y,\mathcal{C}Z).\label{eq:3.10}
\end{eqnarray}
In a similar way, we get
\begin{eqnarray}
\sin^2 \theta g_2((\nabla
F_*)(X,Z_1),F_*(Z_2))&=&-g_1(\mathcal{H}\nabla_{Z_1}\omega \phi
X,Z_2)\nonumber\\
&-&g_1(\mathcal{A}_{Z_1}BZ_2+\mathcal{H}\nabla_{Z_1}\mathcal{C}Z_2,\omega
X).\label{eq:3.11}
\end{eqnarray}
Then proof follows from (\ref{eq:3.10}) and (\ref{eq:3.11}).

\end{document}